\documentclass[12pt]{article}
\font\MSBMx			msbm10
\font\EUFMx			eufm10
   \font\EUFMvii			eufm7
    \font\EUFMv			eufm5
\usepackage{amssymb,latexsym,
amsfonts,amsbsy,color}

\newfam\German
\textfont\German=\EUFMx
\scriptfont\German=\EUFMvii
\scriptscriptfont\German=\EUFMv

\newfam\Blackboadbold
\textfont\Blackboadbold=\MSBMx
\scriptfont\Blackboadbold=\MSBMx

\newcommand{\qed}{\hbox{\rule[-2pt]
{3pt}{6pt}}}

\newtheorem{dfe}{Definition}
[section]

\newtheorem{theo}[dfe]{Theorem}

\newtheorem{pro}[dfe]{Proposition}
\newtheorem{cor}[dfe]{Corollary}

\makeatletter

\@addtoreset{equation}{section}
\makeatother



\begin{title}
{Tight $9$-designs on two concentric spheres}
\end{title}
\author{Eiichi Bannai\quad and \quad
Etsuko Bannai\\
}

\date{}
\begin{document}
\maketitle
\begin{abstract}The main purpose of this paper is to show the nonexistence of tight
Euclidean 9-designs on 2 concentric spheres in $\mathbb R^n$ if $n\geq 3.$
This in turn implies the nonexistence of minimum cubature formulas of
degree 9 (in the sense of Cools and Schmid) for any spherically  symmetric
integrals in $\mathbb R^n$ if $n\geq 3.$
\end{abstract}
\section{Introduction}

The concept of Euclidean $t$-designs $(X, w)$, 
a pair of finite set $X$  in $\mathbb R^n$
and a positive weight function $w$ on $X$, is due to Neumaier-Seidel \cite{N-S}, though
similar concepts have been existed in statistics as rotatable designs \cite{B-H}
and in numerical analysis as cubature formulas for  spherically symmetric
integrals in $\mathbb R^n$ (\cite{C-S,B-H}, etc.). There exist natural Fisher type  lower bounds
(M\"oller's bound) for the size of Euclidean t-designs. Those which  attain
one of such lower bounds are called tight Euclidean 
$t$-designs. These  lower bounds
are basically obtained as functions of $t$, $n$ and the number $p$ of  spheres (whose
centers are at the origin) which meet the finite set 
$X$. We have  been working on
the classification of tight Euclidean $t$-designs, in particular those  with $p=2$ (or $p$
being small). In \cite{B-1} and \cite{B-B-5}, we gave the complete classification of  tight Euclidean
$5$- and $7$-designs on $2$ concentric spheres in $\mathbb R^n$.  (Exactly speaking modulo the existence of tight spherical $4$-designs for $t=5$.)
The main purpose  of this paper
is to show the nonexistence of 
tight Euclidean $9$-designs on $2$ concentric spheres
in $\mathbb R^n$ if $n\geq 3.$

The theory of Euclidean $t$-designs has strong connections with the  theory of
cubature formulas for so called spherically symmetric integrals on  $\mathbb R^n.$
Here, we consider a pair 
$(\Omega, d\rho (\boldsymbol x))$ such that $\Omega$
is a symmetric (or sometimes called radially symmetric)
subset of $\mathbb R^n$ and a spherically symmetric
(or radially symmetric) measure $d\rho (\boldsymbol x)$ 
on $\Omega$.
(Here, a subset $\Omega\subset \mathbb R^n$ is
called spherically symmetric if $\boldsymbol x\in 
\Omega$, then any elements having the same
distance from the origin as $\boldsymbol x$ are 
also in $\Omega$, and $d\rho (\boldsymbol x)$ is 
spherically symmetric if it is invariant under the action
of orthogonal transformations.) A cubature formula $(X, w)$ of degree $t$ for $(\Omega, d\rho (\boldsymbol x))$  is defined as follows.

$X$ is a subset in $\Omega$ containing a finite number of points, $w$ is a positive weight function of $X$, i.e., a map from $X$ to $\mathbb R_{>0}$, and $(X,w)$ satisfies the following condition:
$$\int_\Omega f(\boldsymbol x)d\rho(\boldsymbol x)
=\sum_{\boldsymbol x\in X}w(\boldsymbol x)f(\boldsymbol x)$$
for any polynomials $f(\boldsymbol x)$ of degree at most $t$.

Natural lower bounds of the size $|X|$ of a cubature formula $(X, w)$ of degree $t$ for spherically
symmetric $(\Omega, d\rho (\boldsymbol x))$ are known as M\"oller's lower bounds as follows (\cite{M-1,M-2}).

\begin{enumerate}
\item If $t=2e$, then
$$|X|\geq \dim(\mathcal P_e(\Omega)).$$
\item If $t=2e+1$, then
$$|X|\geq \left\{
\begin{array}{ll} 2\dim(\mathcal P^*_e(\Omega))-1&
\mbox{if $e$ is even and $\boldsymbol 0\in X$,}\\
2\dim(\mathcal P^*_e(\Omega))&\mbox{otherwise.}
\end{array}\right.
 $$ 
 \end{enumerate}
In above $\mathcal P_e(\mathbb R^n)$
is the vector space of polynomials of degree
at most $e$ and $\mathcal P_e(\Omega)=
\{f|_\Omega\mid f\in \mathcal P_e(\mathbb R^n)\}$, and
$\mathcal P_e^*(\mathbb R^n)$
is the vector space of polynomials whose terms are all
of degrees with the same parity as $e$ and
at most $e$. Also $\mathcal P_e^*(\Omega)=
\{f|_\Omega\mid f\in \mathcal P_e^*(\mathbb R^n)\}$.

It is called a minimal cubature formula of degree $t$, if it satisfies a M\"oller's lower bound.
Finding and classifying minimal cubature formulas have been interested  by many
researchers in numerical analysis, and have been studied  considerably
(see \cite{C-S,H-S-1,H-S-2,X}, etc.). As it was pointed out by Cools-Schmid \cite{C-S}, the  problem has a
special feature when $t=4k+1.$ In this case, we can conclude that
(1) $\boldsymbol 0\in X$,
(2) $X$ is on $k+1$ concentric spheres, including $S_1=\{\boldsymbol 0\}.$

Cools-Schmid \cite{C-S} (cf. also \cite{V-C}) gave a complete determination of minimal cubature  formulas
for $n=2$ when $t=4k+1$. The case of $t=5$
for arbitrary $n$ was solved by Hirao-Sawa \cite{H-S-1}
completely, in the effect that the existence of minimal cubature  formula (for
any spherically symmetric $(\Omega, d\rho(\boldsymbol x))$ in 
$\mathbb R^n$ is equivalent to the
existence of tight spherical $4$-design in 
$\mathbb R^n.$ More recently, Hirao-Sawa \cite{H-S-1}
discusses the case of $t=9$ for many specific classical 
$(\Omega, d \rho (\boldsymbol x)).$
As a corollary of our main theorem:  nonexistence of tight
Euclidean $9$-designs on $2$ concentric spheres in 
$\mathbb R^n$ if $n\geq 3,$
we obtain the nonexistence of minimum cubature formulas of
degree $9$ (in the sense of Cools and Schmid) for any spherically  symmetric
integrals in $\mathbb R^n$ if $n\geq 3.$ So, we think that this means a  usefulness of
the concept of Euclidean $t$-design as a master class for all  spherically
symmetric cubature formulas. 
At the end, we add our hope to study the
classification problems of tight Euclidean $t$-designs (for larger $t$) on
2 concentric spheres (or $p$ concentric spheres with small $p$), and  to study
minimal cubature formulas with $t=4k+1$ for $t\geq 13$, 
extending the method used in the present paper.

For more information on spherical designs, Euclidean designs, 
please refer \cite{B-B-1},\cite{B-B-6}, etc..  
Explicit examples of tight $4$-, $5$-, $7$- designs 
on 2 concentric spheres are given in \cite{B-2}, \cite{B-1}, 
\cite{B-B-5}, etc.  
\\

The following is the main theorem of this paper.
\begin{theo}\label{theo:1-1}
Let $(X,w)$ be a tight $9$-design
on $2$ concentric spheres in $\mathbb R^n$ 
of positive radii. 
Let $X=X_1\cup X_2$.
Then the following hold.
\begin{enumerate}
\item $X$ is antipodal.
\item Let $\boldsymbol x\in X_1$, $\boldsymbol y\in X_2$.
Then $\frac{\boldsymbol x\cdot \boldsymbol y}{r_1r_2}$
is a zero of the Gegenbauer polynomial
$Q_{4,n-1}(x)$ of degree $4$. More explicitly,
$Q_{4,n-1}(x)=\frac{n(n+6)}{24}((n+4)(n+2)x^4
-6(n+2)x^2+3)
$ {\rm (Here Gegenbauer polynomial $Q_{l,n-1}(x)$ of degree 
$l$ is normalized so that $Q_{l,n-1}(1)$ is the dimension of
the vector space of homogeneous polynomials of degree $l$.)}.
\item $n=2$ and $(X,w)$ must be similar to the following.

$Y=Y_1\cup Y_2$, $Y_1$ and $Y_2$ are regular $8$-gons given by
\begin{eqnarray}&&Y_1=\left\{r_1(\cos\theta_k,\sin\theta_k)
\mid \theta_k=\frac{2 k\pi}{8}, 0\leq k\leq 7
\right\},
\nonumber\\
&&Y_2=\left\{r_2(\cos\theta_k,\sin\theta_k)
\mid \theta_k=\frac{(2k+1)\pi }{8}, 0\leq k\leq 7
\right\},
\nonumber
\end{eqnarray}
where $r_1$ and $r_2$ are any positive real number 
satisfying $r_1\neq r_2$.
The weight function is defined by $w(\boldsymbol y)=w_1$ on $Y_1$ 
and 
 $w(\boldsymbol y)=\frac{r_1^8}{r_2^8}w_1$ on $Y_2$.
\end{enumerate}
\end{theo}

 It is known that tight Euclidean $(2e+1)$-designs of 
 $\mathbb R^n$
containing the origin exist only when 
$e$ is an even integer and $p=\frac{e}{2}+1$
(see Proposition 2.4.5 in \cite{B-B-H-S}).
Hence Theorem \ref{theo:1-1} implies the followings.
\begin{cor}
Let $(X,w)$ be a tight $9$-design of 
$\mathbb R^n$ containing the origin. Then $n=2$ and
$X$ is supported by $3$ concentric spheres and
$(X\backslash\{\boldsymbol 0\},w)$
is similar to the $9$-design $(Y,w)$ given in Theorem 1.1.
\end{cor}
\begin{cor} If $n\geq 3$, then there is no cubature formula 
of degree $9$ for spherically symmetric subset and measure $(\Omega, d\rho(\boldsymbol x))$ in $\mathbb R^n$.
{\rm(For minimal cubature formulas for $n=2$ see \cite{H-S-2}.) }
\end{cor}

\section{Definition and basic facts on the Euclidean $t$-designs}
We use the following notation.\\
Let $\mathcal P(\mathbb R^n)$ be the vector 
space over real number field $\mathbb R$
consists of all the polynomials in $n$ variables
$x_1,x_2,\ldots,x_n$ with real valued coefficients.
For $f\in \mathcal P(\mathbb R^n)$, $\deg(f)$
denotes the degree of the polynomial $f$.
Let $\mbox{Harm}(\mathbb R^n)$ the subspace of
$\mathcal P(\mathbb R^n)$ consists of all the harmonic
polynomials.
For each nonnegative integer $l$,
let $\mbox{Hom}_l(\mathbb R^n)=\langle f\in \mathcal P(\mathbb R^n)\mid \deg(f)=l\rangle$.
We use the following notation:
\begin{eqnarray}&&\mbox{Harm}_l(\mathbb R^n):=\mbox{Harm}(\mathbb R^n)\cap \mbox{Hom}_l(\mathbb R^n),
\qquad\mathcal P_e(\mathbb R^n):=\oplus_{l=0}^e \mbox{Hom}_l(\mathbb R^n),\nonumber\\
&&\mathcal P_e^*(\mathbb R^n)
:=\oplus_{l=0}^{[\frac{e}{2}]} \mbox{Hom}_{e-2l}(\mathbb R^n),
\quad \mathcal R_{2(p-1)}(\mathbb R^n):=\langle \|\boldsymbol x\|^{2i}\mid 0\leq i\leq p-1\rangle
\subset \mathcal P_{2(p-1)}(\mathbb R^n)
\nonumber
\end{eqnarray}
For a subset $Y\subset \mathbb R^n$,
$\mathcal P(Y)=\{f|_Y\mid f\in \mathcal P(\mathbb R^n)\}$. $\mathcal H(Y)$,  $\mbox{Hom}_l(Y)$,
$\mbox{Harm}_l(Y)$, ...., etc., are defined
in the same way.

Let $(X,w)$ be a weighted finite set in $\mathbb R^n$
whose weight satisfies $w(\boldsymbol x)>0$
for $\boldsymbol x\in X$.
Let $\{r_1,\ r_2,\ldots, r_p\}$ be the set $\{\|\boldsymbol x\| \mid
 \boldsymbol x\in X\}$ of the
length of the vectors in $X$. Where for $\boldsymbol x=(x_1,x_2,\ldots,x_n)$,
$\boldsymbol y=(y_1,y_2,\ldots,y_n)\in \mathbb R^n$, $\boldsymbol x\cdot \boldsymbol y
=\sum_{i=1}^nx_iy_i$ and $\|\boldsymbol x\|
=\sqrt{\boldsymbol x\cdot \boldsymbol x}$.
Let $S_i,\ 1\leq i \leq p,$ be the sphere of radius $r_i$
centered at the origin. 
We say that $X$ is supported by $p$ concentric spheres, or the union of $p$ concentric spheres
$S=S_1\cup S_2\cup\cdots\cup S_p$.

If a finite positive weighted set $(X,w)$ is supported by $p$ concentric spheres, then
$\dim(\mathcal R_{2(p-1)}(X))=p$ holds.
For each $l$, we define an inner product
$\langle -,\ - \rangle_l$ on $\mathcal P_{2(p-1)}(X)$
by 
$\langle f, g\rangle_l
=\sum_{\boldsymbol x\in X}
w(\boldsymbol x)\|\boldsymbol x\|^{2l}
f(\boldsymbol x)g(\boldsymbol x)$.
Then $\langle -,\ -\rangle_l$ is a positive
definite for each $l$.
For each $l$, we define polynomials 
$\{g_{l,j}\mid 0\leq j\leq p-1\}
\subset \mathcal R_{2(p-1)}(\mathbb R^n)$
so that $\{g_{l,j}|_X\mid 0\leq j\leq p-1\}$
is an orthnomal basis of 
$\mathcal R_{2(p-1)}(X)$ with respect 
$\langle -,\ -\rangle_l$.
We define so that $g_{l,j}(\boldsymbol x)$ is a polynomial of degree $2j$ and a linear combination 
of $\{\|\boldsymbol x\|^{2i}\mid 0\leq i\leq j\}$.
We abuse the notation and we identify
 $g_{l,j}(\boldsymbol x)=g_{l,j}(r_\nu)$
 for $\boldsymbol x\in X_\nu\ (1\leq \nu\leq p)$.
 
\begin{dfe}[\cite{N-S}] A weighted finite set $(X,w)$ is a Euclidean
$t$-design if 
$$\sum_{i=1}^p\frac{w(X_i)}{|S_i|}\int_{S_i}
f(\boldsymbol x)d\sigma_i(\boldsymbol x)
=\sum_{\boldsymbol x\in X}w(\boldsymbol x)
f(\boldsymbol x)
$$
holds for any $f\in \mathcal P_t(\mathbb R^n)$.
In above, $w(X_i)=\sum_{\boldsymbol x\in X_i}
w(\boldsymbol x)$, $\int_{S_i}f(\boldsymbol x)d\sigma_i(\boldsymbol x)$ is the usual surface
integral of the sphere $S_i$ of radius $r_i$, 
$|S_i|$ is the surface area of $S_i$. 
\end{dfe}

\begin{theo}[\cite{M-1,M-2,N-S, D-S,B-1,B-B-H-S}, etc]\label{theo:Lb}  Let $X\subset \mathbb R^n$ be a Euclidean $t$-design
supported by a union $S$ of $p$ concentric spheres. Then the following hold.
\begin{enumerate}
\item For $t=2e$,
$$|X|\geq \dim(\mathcal P_e(S)).$$
\item For $t=2e+1$,
$$
|X|\geq\left\{
\begin{array}{ll}
2\dim(\mathcal P^*_e(S))-1
& \mbox{for $e$ even and $\boldsymbol 0\in X$}\\
2\dim(\mathcal P^*_e(S))&
\mbox{otherwise.}
\end{array}
\right.
$$
\end{enumerate}
\end{theo}
\begin{dfe}[Tightness of designs]
If an equality holds in one of the inequalities given
in Theorem \ref{theo:Lb}, then $(X,w)$ is a tight
$t$-design on $p$ concentric spheres in $\mathbb R^n$.
Moreover if $\mathcal P_e(S)=\mathcal P_e(\mathbb R^n)$ holds for $t=2e$, or $\mathcal P_e^*(S)=\mathcal P_e^*(\mathbb R^n)$ holds for $t=2e+1$, then
$(X,w)$ is a tight $t$-design of $\mathbb R^n$. 
\end{dfe}
M\"oller \cite{M-2} proved that
a tight $(2e+1)$-design $(X,w)$ on $p$ concentric spheres is antipodal and the weight function is center symmetric if
$e$ is odd or $e$ is even and $\boldsymbol 0\in X$.
For the case $e$ is even and $\boldsymbol 0\not\in X$  Theorem  2.3.6 in \cite{B-B-H-S} implies
if we assume $p\leq \frac{e}{2}+1$, then $X$
is antipodal and the weight function is center symmetric.
Hence Lemma 1.10 in \cite{B-B-3} and Lemma 1.7 in \cite{B-1} implies
that weight function of a tight $t$-design on $p$ concentric spheres is constant on each $X_i$
for $t=2e$; $t=2e+1$ and $e$ odd;
$t=2e+1$, $e$ even and $\boldsymbol 0\in X$;
$t=2e+1$, $e$ even, $\boldsymbol 0\not\in X$
and $p\leq \frac{e}{2}+1$;

\begin{pro}\label{pro:2-4} Let $(X,w)$ be a positive weighted finite subset in $\mathbb R^n$. 
Assume $\boldsymbol 0\not\in X$ and the weight function is constant on each $X_i$ $(1\leq i\leq p)$. Then the following holds.
$$\sum_{j=0}^{p-1}g_{l,j}(r_\nu)g_{l,j}(r_\mu)=\delta_{\nu,\mu}\frac{1}{|X_\nu|w_\nu r_\nu^{2l}}.$$
\end{pro}
{\bf Proof}\quad Let $M_l$ be 
the $p\times p$ matrix whose $(\nu, j)$ entry
is defined by $\sqrt{|X_\nu|w_\nu} r_\nu^l
g_{l,j}(r_\nu)$ for
$1\leq \nu\leq p,
0\leq j\leq p-1$. Then
\begin{eqnarray}&&({^tM}_l M_l)(j_1,j_2)=\sum_{\nu=1}^p
 M_{\nu,j_1}M_{\nu,j_2}=
 \sum_{\nu=1}^p|X_\nu|w_\nu r_\nu^{2l} 
 g_{l,j_1}(r_\nu)g_{l,j_2}(r_\nu)\nonumber\\
 &&
 =\sum_{\nu=1}^p\sum_{\boldsymbol x\in X_\nu}
 w(\boldsymbol x)\|\boldsymbol x\|^{2l}g_{l,j_1}(r_\nu)g_{l,j_2}(r_\nu)=\sum_{\boldsymbol x\in X^*}
 w(\boldsymbol x)\|\boldsymbol x\|^{2l}g_{l,j_1}(\boldsymbol x)g_{l,j_2}(\boldsymbol x)\nonumber\\
 &&=\delta_{j_1,j_2}
 \end{eqnarray}
 Hence $M_l$ is invertible and $M_l^{-1}={^tM}_l$.
 Hence we have $M_l\ {^tM}_l=I$.
 \begin{eqnarray}
 &&(M_l\ {^tM}_l)(\nu,\mu)=r_\nu^lr_\mu^l\sqrt{|X_\nu||X_\mu|w_\nu w_\mu}\sum_{j=0}^{p-1}g_{l,j}(r_\nu)g_{l,j}(r_\mu)
 =\delta_{\nu,\mu}
 \end{eqnarray}
 Hence we must have
 $$\sum_{j=0}^{p-1}g_{l,j}(r_\nu)g_{l,j}(r_\mu)=
 \delta_{\nu,\mu}\frac{1}{|X_\nu|w_\nu r_\nu^{2l}}$$
 \hfill\qed

\section{Proof of Theorem \ref{theo:1-1} (2)}

Now we prove Theorem 1.1.
Let $(X,w)$ be a tight $9$-design on $2$ 
concentric spheres and $\boldsymbol 0\not\in X$. Let $X=X_1\cup X_2$.
By assumption
$|X|=2\dim(\mathcal P_4^*(S))
=2(\sum_{i=0}^1{n+4-2i-1\choose 4-2i})=\frac{n(n+1)(n^2+5n+18)}{12}$.
Then, as we mentioned in \S 2, $X$ is antipodal and
the weight function is constant on each $X_i$, $i=1,\ 2$. 
Let $w_i=w(\boldsymbol x)$ for $\boldsymbol x\in X_i$.

Let $A(X_i)=\{\frac{\boldsymbol x\cdot  \boldsymbol y}{r_i^2}\mid \boldsymbol x\neq \boldsymbol y\in X_i\}$
 for $i=1,2$.
 Let $A(X_1,X_2)=\{\frac{\boldsymbol x\cdot  \boldsymbol y}{r_1r_2}\mid \boldsymbol x\in X_1,\boldsymbol y\in X_2\}$.
Then
$X_1$ and $X_2$ are spherical $7$-designs and $|A(X_1)|,\ |A(X_2)|\leq 5$ and
 $ |A(X_1,X_2)|\leq 4$. 
 Since $X_1,\ X_2$ are spherical $7$-designs,
 $|X_1|,\ |X_2|\geq \frac{1}{3}(n+2)(n+1)n$.
 We may assume $|X_1|\leq |X_2|$.
 Hence 
 $$\frac{1}{3}(n+2)(n+1)n\leq
 |X_1|\leq \frac{|X|}{2}\leq |X_2|\leq
 |X|-|X_1|\leq \frac{1}{12}n(n+1)(n^2+n+10)
$$
 holds. If $n=2$, then we must have 
 $|X_1|=|X_2|=8$ and
 $X_1$ and $X_2$ are spherical tight $7$-designs.
 We can easily check that for any $A(X_1,X_2)
 =\{\cos(\frac{k\pi}{8})\mid k=1,3,5,7\}
 =\{\frac{\sqrt{2\pm\sqrt{2}}}{2},\
-\frac{\sqrt{2\pm\sqrt{2}}}{2}\}$. Hence $\gamma\in A(X_1,X_2)
$ is a zero of Gegenbauer 
 polynomial $Q_{4,1}(x)=16x^2-16x+2$.\\
 
 In the following we assume $n\geq 3$, then 
$|X_2|\geq \frac{|X|}{2}=\frac{n(n+1)(n^2+5n+18)}{24}> \frac{1}{3}(n+2)(n+1)n$ holds and $X_2$ is not a spherical tight
$7$-design. Hence $X_2$ is a $5$-distance set,
i.e., $|A(X_2)|=5$.
Let $X_i$ be an antipodal half of $X_i^*$ for $i=1,2$.
That is, $X_i=X_i^*\cup (-X_i^*)$, $X_i^*\cap (-X_i^*)=
\emptyset$. 
Then $|A(X_i^*)|\leq 4$ for $i=1,2$, and
$|A(X_1^*,X_2^*)|\leq 4$ hold.

Then equations (3,1) and (3,2) in the proof of Lemma 1.7 in \cite{B-1}  
imply the following equations.

\noindent
$\boldsymbol x\in X_1^*$
\begin{eqnarray}
&&r_1^8g_{4,0}(r_1)^2Q_{4}(1)
+r_1^{4}Q_{2}(1)\sum_{j=0}^1g_{2,j}(r_1)^2
+\sum_{j=0}^1g_{0,j}(r_1)^2
=\frac{1}{w_1}
\end{eqnarray}
$\boldsymbol x\in X_2^*$
\begin{eqnarray}
&&r_2^8g_{4,0}(r_2)^2Q_{4}(1)
+r_2^{4}Q_{2}(1)\sum_{j=0}^1g_{2,j}(r_2)^2
+\sum_{j=0}^1g_{0,j}(r_2)^2
=\frac{1}{w_2}\end{eqnarray}

\noindent
$\boldsymbol x\neq\boldsymbol  y\in X_1^*$
\begin{eqnarray}
&&r_1^8g_{4,0}(r_1)^2
Q_{4}(\frac{(\boldsymbol x,\boldsymbol y)}{r_1^2})
+r_1^4Q_{2}(\frac{(\boldsymbol x,\boldsymbol y)}{r_1^2})
\sum_{j=0}^1g_{2,j}(r_1)^2
+\sum_{j=0}^1g_{0,j}(r_1)^2
=0
\end{eqnarray}

\noindent
$\boldsymbol x\neq \boldsymbol y\in X_2^*$
\begin{eqnarray}
&&r_2^8g_{4,0}(r_2)^2
Q_{4}(\frac{(\boldsymbol x,\boldsymbol y)}{r_2^2})
+r_2^4Q_{2}(\frac{(\boldsymbol x,\boldsymbol y)}{r_2^2})
\sum_{j=0}^1g_{2,j}(r_2)^2
+\sum_{j=0}^1g_{0,j}(r_2)^2=0
\end{eqnarray}

\noindent
$\boldsymbol x\in X_1^*,\boldsymbol  y\in X_2^*$
\begin{eqnarray}
&&r_1^4r_2^4g_{4,0}(r_1)g_{4,0}(r_2)
Q_{4}(\frac{(\boldsymbol x,\boldsymbol y)}{r_1r_2})
+r_1^2r_2^2Q_{2}(\frac{(\boldsymbol x,
\boldsymbol y)}{r_1r_2})
\sum_{j=0}^1g_{2,j}(r_1)g_{2,j}(r_2)\nonumber\\
&&+\sum_{j=0}^1g_{0,j}(r_1)g_{0,j}(r_2)
=0
\end{eqnarray}
In above $g_{l,j}$ are defined for antipodal half $X^*
=X_1^*\cup X_2^*$ of $X$.
Since $X_i^*$ is any antipodal half of $X_i$
for $i=1,\ 2$, 
Proposition \ref{pro:2-4} implies
$$Q_{4,n-1}\left(\frac{\boldsymbol x\cdot\boldsymbol y}{r_1r_2}\right)=0$$
holds for any $\boldsymbol x\in X_1$ and $\boldsymbol y\in X_2$.
\hfill\qed

\begin{pro}\label{pro:3-1} Notation and definition are as given above. $|A(X_1,X_2)|=4$ holds and
$$A(X_1,X_2)=
\left\{\pm\sqrt{\frac{3n+6+\sqrt{6(n+2)(n+1)}}
{(n+4)(n+2)}},\
\pm\sqrt{\frac{3n+6-\sqrt{6(n+2)(n+1)}}
{(n+4)(n+2)}}
\right\}$$
\end{pro}
{\bf Proof}\quad Theorem 1.4 and Theorem 1.5 in \cite{B-B-7} imply that $X$ has the structure of a coherent configuration.
Since $X$ is antipodal
and $0\not\in A(X_1,X_2)$, either $|A(X_1,X_2)|=2$ or 
$|A(X_1,X_2)|=4$ holds.
First assume $|A(X_1,X_2)|=2$. 
Then $A(X_1,X_2)=\{\gamma, -\gamma\}$ with some $\gamma>0$ satisfying $Q_{4,n-1}(\gamma)=0$. Let $\gamma_1=\gamma$
and $\gamma_2=-\gamma$.
Since $X_2$ is a $5$-distance set let
$A(X_2)=\{-1, \pm\beta_2,\pm\beta_4\}$ with real numbers
$\beta_2>\beta_4>0$. Let $\beta_0=1,\ \beta_1=-1$,
$\beta_3=-\beta_2, \beta_5=-\beta_4$.
Then
Proposition 3.2 (1) in \cite{B-B-7} the following hold for any nonnegative integers $l,\ k,\ j$ satisfying $l+k+2j\leq 9$ 
\begin{eqnarray}&&\sum_{u=2}^5
\sum_{v=2}^5
w_2 r_2^{l+k+2j}
Q_{l,n-1}(\beta_u)Q_{k,n-1}(\beta_v)
p_{\beta_u,\beta_v}^{\beta_0}
\nonumber\\
&&
+
\sum_{u=1}^2\sum_{v=1}^2 w_1 r_1^{l+k+2j}
Q_{l,n-1}(\gamma_u)Q_{k,n-1}(\gamma_v)
p_{\gamma_u,\gamma_v}^{\beta_0}
\nonumber\\
&&
=\delta_{l,k}Q_{l,n-1}(1)\sum_{\nu=1}^2
N_\nu w_\nu r_\nu^{2l+2j}\nonumber\\
&&
-w_2r_2^{l+k+2j}((-1)^{l+k}+1)
Q_{l,n-1}(1)Q_{k,n-1}(1),
\end{eqnarray}
$N_\nu=|X_\nu|$ for $\nu=1,\ 2$ and
$p_{\beta_u,\beta_v}^{\beta_0}$, $p_{\gamma_u,\gamma_v}^{\beta_0}$ denotes the corresponding intersection numbers.
 Since $Q_{4,n-1}(\gamma)=Q_{4,n-1}(-\gamma)=0$,
$p_{\beta_u,\beta_v}^{\alpha_0}=0$, for any $2\leq u\neq v\leq 5$, and
$p_{\gamma_u,\gamma_v}^{\alpha_0}=0$, for any $1\leq u\neq v\leq 2$, $p_{\gamma_1,\gamma_1}^{\beta_0}=p_{\gamma_2,\gamma_2}^{\beta_0}=\frac{|X_1|}{2}$,
$p_{\beta_3,\beta_3}^{\beta_0}=p_{\beta_2,\beta_2}^{\alpha_0}$, $p_{\beta_5,\beta_5}^{\beta_0}=p_{\beta_4,\beta_4}^{\alpha_0}$,
equations for
$(l,k,j)=(0,0,0), (1,0,0), (1,1,0), (2,1,1)$ imply
$$p_{\beta_2,\beta_2}^{\beta_0}
= \frac{
-w_2r_2^2(n(N_2-2)\beta_4^2-N_2+2n)
-N_1w_1r_1^2(-1+n\gamma_1^2)}
{2nw_2r_2^2(\beta_2^2-\beta_4^2)}$$
and
$$p_{\beta_4,\beta_4}^{\beta_0}=
\frac{w_2r_2^2
(n(N_2-2)\beta_2^2-N_2+2n)+N_1w_1r_1^2(-1+n\gamma_1^2)}
{2nw_2r_2^2(\beta_2^2-\beta_4^2)}.
$$
Then equation for
$(l,k,j)=(1,1,1)$ implies 
$$(r_1^2-r_2^2)(-1+n\gamma_1^2)r_2^2w_1N_1n=0.$$
Since $\gamma_1$ is a zero of $Q_{4,n-1}(x)$, this is a contradiction.
\hfill\qed\\

Since $n\geq 3$, we have 
$|X_2|\geq \frac{1}{2}|X|
=\frac{1}{24}n(n+1)(n^2+5n+18)>\frac{1}{3}(n+2)(n+1)n$.
We divide the proof of Theorem \ref{theo:1-1} 
into two cases I and II. 
In Case I, we assume $X_1$ is not a tight spherical $7$-design, i.e. $|X_1|>\frac{1}{3}(n+2)(n+1)n$, and
in Case II, we assume $X_1$ is a tight spherical $7$-design,
i.e. $|X_1|=\frac{1}{3}(n+2)(n+1)n$.\\

\noindent
{\bf Case I, $|X_2|\geq |X_1|>\frac{1}{3}(n+2)(n+1)n$}\\
In this case both $X_1$ and $X_2$ are
antipodal spherical $7$-designs and $5$-distance 
sets.
\begin{eqnarray}&&A(X_1)=\{\alpha_1,\ \alpha_2,\ \alpha_3,\ \alpha_4,\ \alpha_5\},\quad 
\alpha_0=1,\ \alpha_1=-1,\ \alpha_3=-\alpha_2,\ \alpha_5=-\alpha_4,\nonumber\\
&&A(X_2)=\{\beta_1,\ \beta_2,\ \beta_3,\ \beta_4,\ \beta_5\},\quad 
\beta_0=1,\ \beta_1=-1,\ \beta_3=-\beta_2,\ \beta_5=-\beta_4,\nonumber\\
&&A(X_1,X_2)=\{\gamma_1,\ \gamma_2,\ \gamma_3,\ \gamma_4\},
\end{eqnarray}
where $\gamma_1=\sqrt{\frac{3n+6+\sqrt{6(n+2)(n+1)}}
{(n+4)(n+2)}},\ \gamma_2=-\gamma_1$,
$\gamma_3=\sqrt{\frac{3n+6-\sqrt{6(n+2)(n+1)}}
{(n+4)(n+2)}},\ \gamma_4=-\gamma_3$. We may assume $\alpha_2>\alpha_4>0$,
 $\beta_2>\beta_4>0$.
Then Proposition 9.1 and Theorem 9.2 in \cite{B-B-5} imply
the followings (see also \cite{B-B-2,B-B-4}).\\
$\bullet$ $X_i^*$ $(1\leq i\leq 2)$ has the structure of
a strongly regular graphs.\\
$\bullet$ $\frac{1-\alpha_2^2}{\alpha_2^2-\alpha_4^2}$
and $\frac{1-\beta_2^2}{\beta_2^2-\beta_4^2}$
are integers.\\
$\bullet$ $\alpha_2,\ \alpha_3,\ \alpha_4,\ \alpha_5$ 
are the zeros of the following polynomial $a(x)$.

$a(x)=(n+4)(n+2)(N_1-n^2-n)x^4+(n+2)( n^3+6n^2+5n-6N_1)x^2+3N_1-n^3-3n^2-2n.
$\\
$\bullet$ $\beta_2,\ \beta_3,\ \beta_4,\ \beta_5$ 
are the zeros of the following polynomial $b(x)$.

$b(x)=(n+4)(n+2)(N_2-n^2-n)x^4+(n+2)( n^3+6n^2+5n-6N_2)x^2+3N_2-n^3-3n^2-2n.
$\\
$\bullet$ $n\geq 4$ and $\alpha_i,$ and $\beta_i$, $i=2,3,4$, are rational numbers.\\
In above $N_i=|X_i|$ for $i=1,\ 2$.

Hence we obtain
\begin{eqnarray}&&\alpha_2^2=\frac{
(n+2)\bigg(6N_1-n(n+1)(n+5)\bigg)+
\sqrt{(n+1)(n+2)D_1}
}{2(n+4)(n+2)(N_1-n^2-n)}
\\
&&\alpha_4^2=\frac{
(n+2)\bigg(6N_1-n(n+1)(n+5)\bigg)-
\sqrt{(n+1)(n+2)D_1}
}{2(n+4)(n+2)(N_1-n^2-n)}
\\
&&\beta_2^2=\frac{
(n+2)\bigg(6N_2-n(n+1)(n+5)\bigg)+
\sqrt{(n+1)(n+2)D_2}
}{2(n+4)(n+2)(N_2-n^2-n)}\label{equ:3-10}
\\
&&\beta_4^2=\frac{
(n+2)\bigg(6N_2-n(n+1)(n+5)\bigg)-
\sqrt{(n+1)(n+2)D_2}}{2(n+4)(n+2)(N_2-n^2-n)}
\label{equ:3-11}
\end{eqnarray}
where $D_1=
n^2(n+1)(n+2)(n+3)^2-8n(n+1)(n+5)N_1+24N_1^2$,
$D_2=
n^2(n+1)(n+2)(n+3)^2-8n(n+1)(n+5)N_1+24N_2^2$
$N_i=|X_i|$ $(1\leq i\leq 2)$.

Next proposition is very important.
\begin{pro} Notation and definition are as given above.
Assume $n\geq 3$, then $\sqrt{6(n+1)(n+2)}$
is an integer and
$\gamma_i^2$ $(1 \leq i\leq 4)$ are rational numbers. 
\end{pro}
{\bf Proof}\quad Theorem 1.4 and Theorem 1.5 in \cite{B-B-7} imply 
that $X$ has the structure of a coherent configuration. 
Let $\boldsymbol x\in X_1$ and $p_{\gamma_i,\gamma_i}^{\alpha_0}
=|\{\boldsymbol z\in X_2\mid 
\frac{\boldsymbol x\cdot \boldsymbol z}{r_1r_2}=\gamma_i\}|$.
Using the equations given in Proposition 3.2 (1) in \cite{B-B-7} the following hold for any nonnegative integers $l,\ k,\ j$ satisfying $l+k+2j\leq 9$ 
\begin{eqnarray}&&\sum_{u=2}^5
\sum_{v=2}^5
w_1r_1^{l+k+2j}
Q_{l,n-1}(\alpha_u)Q_{k,n-1}(\alpha_v)
p_{\alpha_u,\alpha_v}^{\alpha_0}
\nonumber\\
&&
+
\sum_{u=1}^4\sum_{v=1}^4w_2 r_2^{l+k+2j}
Q_{l,n-1}(\gamma_u)Q_{k,n-1}(\gamma_v)
p_{\gamma_u,\gamma_v}^{\alpha_0}
\nonumber\\
&&
=\delta_{l,k}Q_{l,n-1}(1)\sum_{\nu=1}^2
N_\nu w_\nu r_\nu^{2l+2j}\nonumber\\
&&
-w_1r_1^{l+k+2j}((-1)^{l+k}+1)
Q_{l,n-1}(1)Q_{k,n-1}(1)
\end{eqnarray}
Since $p_{\alpha_1,\alpha_1}^{\alpha_0}=1$,
$p_{\alpha_i,\alpha_j}^{\alpha_0}=0$ for any $1\leq i\neq j\leq 5$, and
$p_{\gamma_i,\gamma_j}^{\alpha_0}=0$ for any $1\leq i\neq j\leq 4$, 
we have the followings.
\begin{eqnarray}&&
p_{\gamma_1,\gamma_1}^{\alpha_0}
=p_{\gamma_2,\gamma_2}^{\alpha_0}
= \frac{N_2(1-n\gamma_3^2)
}{2n(\gamma_1^2-\gamma_3^2)},\nonumber\\
&&p_{\gamma_3,\gamma_3}^{\alpha_0}
=p_{\gamma_4,\gamma_4}^{\alpha_0}
= \frac{N_2(n\gamma_1^2-1)}
{2n(\gamma_1^2-\gamma_3^2)}
\end{eqnarray}
Then 
$p_{\gamma_1,\gamma_1}^{\alpha_0}=
\frac{\big(3n^2+3n-(n-2)\sqrt{6(n+1)(n+2)}\big)N_2}{12n(n+1)}$.
Hence $\sqrt{6(n+1)(n+2)}$ is an integer.
This completes the proof.
\hfill\qed\\

Next, we express $\frac{1-\alpha_2^2}{\alpha_2^2-\alpha_4^2}$ and $\frac{1-\beta_2^2}{\beta_2^2-\beta_4^2}$ interms of $n$ and $N_1$, $N_2$.
We have
\begin{eqnarray}
&&\frac{1-\alpha_2^2}{\alpha_2^2-\alpha_4^2}
=-\frac{1}{2}+F(n,N_1),\\
&&\frac{1-\beta_2^2}{\beta_2^2-\beta_4^2}=
-\frac{1}{2}+F(n,N_2),
\end{eqnarray}
where 
\begin{eqnarray}
&&F(n,x)=\nonumber\\
&&\frac{
(2x-n^2-3n)\sqrt{(n+1)(n+2)
\bigg(n^2(n+1)(n+2)(n+3)^2-8n(n+1)(n+5)x+24x^2\bigg)}}
{2\bigg(n^2(n+1)(n+2)(n+3)^2-8n(n+1)(n+5)x+24x^2\bigg)
}\nonumber\\
&&
\end{eqnarray}
We have $\frac{(n+2)(n+1)n}{3}<N_1\leq\frac{1}{24}n(n+1)(n^2+5n+18)
\leq N_2\leq \frac{1}{12}n(n+1)(n^2+n+10)$.
Since 
\begin{eqnarray}
&&F(n,x)=\frac{
(1-\frac{n^2+3n}{2x})}{(
\frac{n^6+9n^5+29n^4+39n^3+18n^2}{2x^2}
-\frac{4n(n^2+6n+5)}{x}+12)
}\times\nonumber\\
&&\nonumber\\
&&\sqrt{6(n+2)(n+1)
(\frac{n^6+9n^5+29n^4+39n^3+18n^2}{24x^2}
-\frac{n(n^2+6n+5)}{3x}+1)},\nonumber
\end{eqnarray}
we can observe that for $x>\frac{1}{24}n(n+1)(n^2+5n+18)$, $F(n,x)\approx\frac{\sqrt{6(n+2)(n+1)}}{12}$. More precisely we have the followings. 
\begin{eqnarray}&&\frac{\partial F(n,x)}{\partial x}=
\frac{(n-1)(n+4)(n+2)(n+1)(n^3+4n^2+3n-4x)n}
{\sqrt{(n+2)(n+1)\bigg(n^2(n+2)(n+1)(n+3)^2-8n(n+5)(n+1)x+24x^2)\bigg)^3}}
\nonumber\\
&&\end{eqnarray}
Hence $ F(n,x)$ decreases for $x\geq \frac{1}{4}n(n+1)(n+3)$.

\begin{eqnarray}&&F(n,\frac{1}{12}n(n+1)(n^2+n+10)
)=\frac{\sqrt{6}(n^2+3n+8)}{12\sqrt{n^2-n+4}}
\nonumber\\
&&> \frac{\sqrt{6(n+1)(n+2)}}{12}
\end{eqnarray}
\begin{eqnarray}&&F(n,\frac{1}{24}n(n+1)(n^2+5n+18))=\frac{\sqrt{6(n+2)}(n^2+7n+18)}
{12\sqrt{n^3+5n^2+16n+36}}\nonumber\\
&&<1+\frac{\sqrt{6(n+1)(n+2)}}{12}
\end{eqnarray}
Hence 
$$-\frac{1}2+\frac{\sqrt{6(n+1)(n+2)}}{12}<-\frac{1}2+F(n,N_2)<\frac{1}2+\frac{\sqrt{6(n+1)(n+2)}}{12}
$$
holds.
Since $\sqrt{6(n+1)(n+2)}$ is an integer,
$\sqrt{6(n+1)(n+2)}=\sqrt{6^2k^2}=6k$ with an
integer $k>0$.
Hence 
$$\frac{k-1}{2}<-\frac{1}2+F(n,N_2)<\frac{k+1}{2}
$$
If $k$ is an odd integer, then $-\frac{1}{2}+F(n,N_2)$ cannot be an integer.
Hence $k$ must be an even integer
and we must have
\begin{eqnarray}&&-\frac{1}{2}+F(n,N_2)=\frac{k}{2}=\frac{\sqrt{6(n+2)(n+1)}}{12}.
\label{equ:3-19}
\end{eqnarray}
It is known $n=23$, $2399$, $235223$ satisfy this condition.
Otherwise $n>300000$.
The equation (\ref{equ:3-19}) implies
\begin{eqnarray}&&N_2=
\frac{n}{36(2n^2+6n+1)}\times\nonumber\\
&&\bigg\{9(n+3)(n+1)(n^2+6n+2)
+(n-1)(n+4)(n+2)(n+1)\sqrt{6(n+1)(n+2)}
\nonumber\\
&&
+\varepsilon(n-1)\bigg(\sqrt{6}(n^2+3n-1)+3\sqrt{(n+2)(n+1)}\bigg)\sqrt{(n+4)(n+1)}\times
\nonumber\\
&&
\sqrt{(n+5)(n+1)-\sqrt{6(n+2)(n+1)}}
\bigg\}\nonumber\\
&&
\end{eqnarray}
where $\varepsilon=1$ or $-1$.
If $\varepsilon=-1$,
then
we have 
$$N_2<\frac{1}{24}n(n+1)(n^2+5n+18).$$
This contradicts the assumption. Hence we must have
$\varepsilon=1$.
Then we must have
\begin{eqnarray}&&N_1=\frac{n}{36(2n^2+6n+1)}\times\nonumber\\
&&
\bigg\{3n(n+1)(2n^3+13n^2+40n+53)
-(n-1)(n+4)(n+2)(n+1)\sqrt{6(n+1)(n+2)}\nonumber\\
&&-(n-1)\bigg(\sqrt{6}(n^2+3n-1)+3\sqrt{(n+2)(n+1)}\bigg)
\sqrt{(n+4)(n+1)}\times\nonumber\\
&&\sqrt{(n+5)(n+1)-\sqrt{6(n+1)(n+2)}}
\bigg\}\nonumber\\
&&
\end{eqnarray}
Since $n=23,\ 2399$, and $235223$ do not give
integral value for $N_2$, we must have
$n>300000$.
Solve 
$-\frac{1}{2}+F(n,x)=\frac{\sqrt{6(n+2)(n+1)}}{12}+2$
for $x$, then
we must have $x=K_\varepsilon$ given below.
\begin{eqnarray}&&K_\varepsilon=
\frac{n}{60(6n^2+18n-213
)}
\times\nonumber\\
&&\bigg\{45(n+1)(n^3+9n^2-28n-234)
+(n-1)(n+4)(n+2)(n+1)\sqrt{6(n+2)(n+1)}
\nonumber\\
&&+\varepsilon(n-1)\bigg(\sqrt{6}(n^2+3n-73)+15\sqrt{(n+2)(n+1)}\bigg)
\times\nonumber\\
&&\sqrt{n^2+6n-67-5\sqrt{6(n+2)(n+1)}}\bigg\}
\end{eqnarray}
where $\varepsilon=\pm $. Now we may assume 
$n>300000$.
Then we have
\begin{eqnarray}
&&K_+>\frac{n}{60(6n^2+18n-213)}\times(n-1)(n+4)(n+2)(n+1)\sqrt{6(n+2)(n+1)}\nonumber\\
&&>\frac{\sqrt{6}n^5(n-1)}{60(6n^2+18n-213)}>\frac{n(n+1)(n+3)}{4}.
\end{eqnarray}
Next compare $K_+$ and $N_1$. 
\begin{eqnarray}
&&N_1-K_+=
\frac{n(n-1)}{180(2n^2+6n+1)(2n^2+6n-71)}\times\nonumber\\
&&
\bigg\{15(n+2)(n+1)(4n^4+28n^3-76n^2-442n-351)
\nonumber\\
&&
-6(n+4)(n+2)(n+1)(2n^2+6n-59)
\sqrt{6(n+2)(n+1)}\nonumber\\
&&
-(2n^2+6n+1)\bigg(\sqrt{6} (n^2+3n-73)
+15\sqrt{(n+2)(n+1)}\bigg)
\times\nonumber\\
&&\sqrt{(n+4)(n+1)}\sqrt{n^2+6n-67
-5\sqrt{6(n+2)(n+1)}}\nonumber\\
&&
-5(2n^2+6n-71)
\bigg(\sqrt{6}(n^2+3n-1)+3\sqrt{(n+2)(n+1)}\bigg)
\times\nonumber\\
&&
\sqrt{(n+4)(n+1)}\sqrt{n^2+6n+5-\sqrt{6(n+2)(n+1)}}\bigg\}
\end{eqnarray}
The order of the  the formula in $\bigg\{ \cdots \bigg\}$ in above 
equals $2(30-11\sqrt{6})n^6$.
Hence $N_1>K_+$ holds for any $n$ sufficiently large, in particular for $n>300000$. 
This means 
$$-\frac{1}{2}+F(n,N_1)<\frac{\sqrt{6(n+2)(n+1}}{12}+2$$
holds for any $n$  sufficiently large. 
Since $N_2>N_1$, we must have 
$\frac{\sqrt{6(n+2)(n+1}}{12}=-\frac{1}{2}+F(n,N_2)
<-\frac{1}{2}+F(n,N_1)$.
Hence we must have 
$-\frac{1}{2}+F(n,N_1)=\frac{\sqrt{6(n+2)(n+1}}{12}+1$.
Next solve for
$F(n,x)=\frac{\sqrt{6(n+2)(n+1)}}{12}+1$
then we have $x=G_\varepsilon$ given below.
\begin{eqnarray}&&
G_\varepsilon=
\frac{n}{6n^2+18n-69}
\times\nonumber\\
&&\bigg\{27(n+1)(n^3+9n^2+4n-74)
+(n-1)(n+4)(n+2)(n+1)\sqrt{6(n+2)(n+1)}
\nonumber\\
&&+\varepsilon(n-1)\bigg(\sqrt{6}(n^2+3n-25)
+9\sqrt{(n+2)(n+1)}\bigg)\times\nonumber\\
&&
\sqrt{(n+4)(n+1)(n^2+6n-19-3\sqrt{6(n+2)(n+1)})})
\bigg\}
\end{eqnarray}
where $\varepsilon=\pm $.
Compare $N_1$ and $G_+$.
\begin{eqnarray}
&&G_+-N_1=\frac{n(n-1)}{108(2n^2+6n+1)(2n^2+6n-23)}\times\nonumber\\
&&\bigg\{-9(n+2)(n+1)(4n^4+28n^3+20n^2-106n-111)
\nonumber\\
&&+4(n+4)(n+2)(n+1)(2n^2+6n-17)\sqrt{6(n+2)(n+1)}
\nonumber\\
&&+(2n^2+6n+1)\bigg(\sqrt{6}(n^2+3n-25)
+9\sqrt{(n+2)(n+1)}\bigg)\times
\nonumber\\
&&
\sqrt{(n+4)(n+1)(n^2+6n-19-3\sqrt{6(n+2)(n+1)})}
\nonumber\\
&&+3(2n^2+6n-23)\bigg(\sqrt{6}(n^2+3n-1)+3\sqrt{(n+2)(n+1)}\bigg)
\times\nonumber\\
&&\sqrt{(n+4)(n+1)\bigg(n^2+6n+5-\sqrt{6(n+2)(n+1)}\bigg)
}\bigg\}
\end{eqnarray}
The order of the formula in $\bigg\{\cdots\bigg\}$
given above equals $4(4\sqrt{6}-9)n^6$. 
Hence $G_+>N_1$ holds for any $n$ sufficiently large, in particular $n>300000$.
Since $F(n,x)$ decreases for $x\geq \frac{(n+3)(n+1)n}{4}$, we have
$$N_2>G_+>N_1>K_+>\frac{(n+3)(n+1)n}{4}.$$
Hence we must have
\begin{eqnarray}&&\frac{\sqrt{6(n+2)(n+1)}}{12}=-\frac{1}{2}+F(n,N_2)
<-\frac{1}{2}+F(n,G_+)=\frac{\sqrt{6(n+2)(n+1)}}{12}+1
\nonumber\\
&&<-\frac{1}{2}+F(n,N_1)<
\frac{\sqrt{6(n+2)(n+1)}}{12}+2.
\end{eqnarray}
Hence, $-\frac{1}{2}+F(n,N_1)$ cannot be an integer for any sufficiently large $n$, in particular for $n>300000$.\\

\noindent
{\bf Case II, $|X_2|> |X_1|=\frac{1}{3}(n+2)(n+1)n$}\\
In this case we must have $|X_2|=\frac{1}{12}n(n+1)(n^2+n+10)$. Since $X_1$ is a tight spherical $7$-design,
$X_1$ is a $4$-distance set. On the other hand $X_2$ is a $5$-distance set.
It is known that $A(X_1)
=\left\{0, -1, \pm\sqrt{\frac{3}{n+4}}\right\}$,
$\sqrt{\frac{n+4}{3}}$ is an integer.
Let $\alpha_1=-1,\ \alpha_2=0,\
\alpha_3=\sqrt{\frac{3}{n+4}},\ 
\alpha_4=-\sqrt{\frac{3}{n+4}}$ 
and $\alpha_0=1$.
By Proposition \ref{pro:3-1}, we have
$\gamma_1=\frac{\sqrt{3n+6+\sqrt{6(n+2)(n+1)}}}{\sqrt{(n+4)(n+2)}}
$,
$\gamma_3=\frac{\sqrt{3n+6-\sqrt{6(n+2)(n+1)}}}{\sqrt{(n+4)(n+2)}}
$.
Proposition 9.1 and Theorem 9.2 in \cite{B-B-5}
imply that (\ref{equ:3-10}) and (\ref{equ:3-11})
also hold in this case. 
Since $N_2=|X_2|=\frac{1}{12}n(n+1)(n^2+n+10)$,
we obtain
$\beta_2=\frac{
\sqrt{(n+4)(n+2)(3n+\sqrt{6n^2-6n+24})}}{(n+4)(n+2)}$
\quad and\quad 
$\beta_4=\frac{
\sqrt{(n+4)(n+2)(3n-\sqrt{6n^2-6n+24})}}{(n+4)(n+2)}$.
Hence we have
$\frac{1-\beta_2^2}{\beta_2^2-\beta_4^2}=
-\frac{1}{2}+\frac{n^2+3n+8}{2\sqrt{6n^2-6n+24}}$.
Therefore
$$-\frac{1}{2}+\frac{n^2+3n+8}{ 2\sqrt{6n^2-6n+24}}$$
is an integer. Then $24(\frac{n^2+3n+8}{ 2\sqrt{6n^2-6n+24}})^2$ must be an integer. Since
$$24\left(\frac{n^2+3n+8}{ 2\sqrt{6n^2-6n+24}}\right)^2=\frac{(n^2+3n+8)^2}{n^2-n+4}=n^2+7n+28
+\frac{48(n-1)}{n^2-n+4}
$$
there is no integer $n$ satisfying the condition. 
This implies that for $n\geq 3$, there is no tight
$9$-design on two concentric spheres satisfying
$N_1=\frac{(n+2)(n+1)n}{3}$.

\end{document}